\documentclass{article}
\usepackage{amsmath, amssymb, latexsym, amsfonts, mathrsfs}
\def\Pr{{\bf P}}
\def\LL{{\cal L}}
\def\sr{{\sf r}}
\def\sC{{\sf C}}
\def\CC{{\cal C}}
\def\CCC{{\mathscr C}}
\def\sM{{\sf M}}
\def\MM{{\cal M}}
\def\MMM{{\mathscr M}}
\def\kbar{{\overline{k}}}

\def\vph{\varphi\0}
\def\ra{\rightarrow}

\def\0{^{\phantom0}}
\def\9{_{\phantom9}}

\def\GV{Gilbert-\mbox{\kern-.16em}Varshamov}
\def\vladut{Vl\u{a}du\c{t}}
\def\DV{Drinfeld-\kern-.16em\vladut}

\begin{document}

\centerline{\large Still better nonlinear codes from modular curves}
\vspace*{2ex}
\centerline{Noam D.\ Elkies}
\centerline{July 2003}

\vspace*{5ex}

\begin{quote}

{\bf Abstract.}
We give a new construction of nonlinear error-correcting codes
over suitable finite fields~$k$\/ from the geometry of modular curves
with many rational points over~$k$,
combining two recent improvements on Goppa's construction.
The resulting codes are asymptotically the best currently known.
\end{quote}

\vspace*{3ex}

{\large\bf 1. Introduction.}

{\bf 1.1 Review of Goppa's construction.}
Fix a finite field~$k$ of $q=p^\alpha$ elements.
Let $C$\/ be a (projective, smooth, irreducible)
algebraic curve of genus~$g$ defined over~$k$,
with $N$\/ rational points.  It is known that
\begin{equation}
\label{eq:DV}
N < (q^{1/2}-1+o(1)) g
\end{equation}
as $g\ra\infty$ (\DV~\cite{DV}).
We say a curve of genus $g\ra\infty$ is ``asymptotically optimal''
if it has at least $(q^{1/2}-1-o(1)) \, g$ rational points over~$k$.
If $\alpha$ is even (that is, if $q_0 := \sqrt q$ is an integer),
then modular curves of various flavors ---
classical (elliptic), Shimura, or Drinfeld --- attain
\begin{equation}
\label{eq:modular}
N \geq (q_0-1)(g-1) = (q^{1/2}-1-o(1)) \, g
\end{equation}
\cite{Ihara, TVZ}, and are thus asymptotically optimal.

Let $D$\/ be a divisor on~$C$\/ of degree $<N$.
Goppa (\cite{Goppa}, see also \cite{TV})
regards the space $\LL(D)$ of sections of~$D$\/
as a linear code in $k^N$,
whose dimension~$\sr$ and minimal distance~$d$\/ satisfy
$\sr \geq \deg(D)-g+1$ (by the Riemann-Roch theorem)
and $d\geq N-\deg(D)$
(because a nonzero section of~$D$\/ has at most $\deg(D)$ zeros).
Thus the transmission rate $R=\sr/N$\/
and the error-detection rate $\delta=d/N$\/ of Goppa's codes
are related by
\begin{equation}
\label{eq:Goppa}
R + \delta > 1 - \frac{g}{N} .
\end{equation}
This lower bound improves as $N/g$ increases.
For $q=q_0^2$, we may take $C$\/ asymptotically optimal, and find
\begin{equation}
\label{eq:Goppamax}
R + \delta > 1 - \frac{1}{q_0-1} - o(1)
\end{equation}
for an infinite family of linear codes over~$k$,
which is the best that can be obtained from~(\ref{eq:Goppa}).

Let us say that $(R_0,\delta_0)$ is {\em asymptotically feasible}
if $R_0,\delta_0$ are positive
and there exist arbitrarily long codes over~$k$, linear or not,
with\footnote{
  As usual, the rate of a nonlinear code $C \subset k^n$
  is defined by $R = N^{-1}\log_q(\#C)$, which equals $\sr/N$\/
  when $C$\/ is linear.
  }
$R>R_0$ and $\delta>\delta_0$.
Then Goppa's construction yields the asymptotic feasibility
of $(R_0,\delta_0)$ for any positive $R_0,\delta_0$
such that $R_0 + \delta_0 < 1 - (1/(q_0-1))$.
This is true because $\deg(D)$ is an arbitrary integer in $(0,N)$,
and nontrivial (in the sense that such $(R_0,\delta_0)$ exist)
once $q_0>2$.

By comparison, a random code or random linear code
of length $N\ra\infty$ and rate $R$\/
has error-detection rate at least $\delta - o(1)$ with probability
$1-o(1)$ provided $\delta < (q-1)/q$ and $R + H_q(\delta) < 1$,
where $H_q$ is the normalized entropy function
\begin{eqnarray}
\label{eq:entropy}
H_q(\delta) \!\! & := & \!\! \delta \log_q(q-1)
- \delta \log_q \delta - (1-\delta) \log_q (1-\delta)
\\
\!\! & = & \!\!
\delta \log_q \left( (q-1) \frac{1-\delta}{\delta} \right)
- \log_q(1-\delta)
\label{eq:entropy1}
\\
\nonumber
& \biggl[ = \, & \!\!
\lim_{N_\ra\infty} N^{-1} \log_q \binom{N}{\delta N} \biggr].
\end{eqnarray}
Therefore if $0 < \delta_0 < (q-1)/q$ and $0 < R_0 < 1 - H_q(\delta_0)$
then $(R_0,\delta_0)$ is asymptotically feasible.
This is the {\em \GV\ bound}.

Once \hbox{$q_0\geq7$}, Goppa's construction yields
asymptotically feasible $(R_0,\delta_0)$ beyond the \GV\ bound.
This was the first construction to improve on \GV\ in this sense,
and it remains the only such construction that yields linear codes.

{\bf 1.2 Beyond Goppa.} 
For about 20 years Goppa's technique remained the best construction
of codes over an alphabet of size $q_0^2 \geq 49$ beyond the
\GV\ bound.  Refinements concerned only algorithmic improvements,
for exhibiting suitable curves~$C$\/ (see
\cite{E:tower} for classical and Shimura curves,
\cite{E:shimura} for further Shimura curves, and
\cite{GS:drinfeld, GS:more_drinfeld, E:drinfeld, TV}
for Drinfeld modular curves)
and using the resulting codes for error-resistant communication
(polynomial-time encoding and decoding, see~\cite{GuSu,ShWa}).
In \cite{E:stoc},\footnote{
  While \cite{E:stoc} was the first publication,
  we obtained these results in the mid-1990's and included them
  in several conference and seminar talks starting in~1996.
  }
we used rational functions on~$C$\/
to construct algebraic-geometry codes over the $(q+1)$-letter alphabet
$\Pr^1(k) = k \cup \{ \infty \}$, and gave asymptotic estimates
on their parameters $R,\delta$.  We argued in~\cite{E:stoc}
that these codes improve on Goppa's in a range of parameters
that includes all the Goppa codes that improve on~\GV;
but our comparison was necessarily indirect
due to the different alphabet sizes.

Even more recently, Xing~\cite{Xing} gave a new construction
of nonlinear algebraic-geometry codes over~$k$.
Like Goppa, Xing uses sections of line bundles,
but he cleverly exploits the sections' derivatives
to find codes with better asymptotic parameters than Goppa's.
The Xing codes have
\begin{equation}
\label{eq:Xing}
R + \delta > 1 - \frac{1}{q_0-1}
  + \sum_{i=2}^\infty \log_q \left( 1 + \frac{q-1}{q^{2i}} \right)
  - o(1),
\end{equation}
which improves on (\ref{eq:Goppamax}) by
\begin{equation}
\label{eq:c_q}
c_q := \sum_{i=2}^\infty \log_q \left( 1 + \frac{q-1}{q^{2i}} \right)
= \frac1{\log q}(q^{-3} - q^{-4} + O(q^{-5})).
\end{equation}
In particular, this is the first construction
of algebraic-geometry codes over a \hbox{$4$-letter} alphabet
with $R,\delta$ both bounded away from~zero.
(Our codes over~$\Pr^1(k)$ attained this for an alphabet of $5$~letters.)
Xing's construction does not require that $q=q_0^2$,
and yields an improvement of~$c_q$ over the Goppa bound (\ref{eq:Goppa})
for all~$q$.

In this paper we obtain a further improvement, at least for $q=q_0^2$,
by applying Xing's technique to our codes of~\cite{E:stoc}.
While those codes used an alphabet of $q+1$ letters,
our new codes $\CCC_D(h)$ use the \hbox{$q$-letter} alphabet~$k$,
and can thus be compared directly with Goppa's and Xing's codes.
We find that when $q=q_0^2$ our $\CCC_D(h)$ have parameters
that improve on Xing's, replacing the sum (\ref{eq:c_q}) by
\begin{equation}
\label{eq:c_q+}
\log_q \left( 1 + \frac1{q^3} \right)
= \frac1{\log q}(q^{-3} - O(q^{-6})).
\end{equation}
We conclude that $(R_0,\delta_0)$ is asymptotically feasible
for any positive $R_0,\delta_0$ such that
\begin{equation}
\label{eq:feasible}
R_0 + \delta_0 < 1 - \frac1{q_0-1}
+ \log_q \left( 1 + \frac1{q^3} \right).
\end{equation}
This also gives further support to our claim
that the codes of~\cite{E:stoc} improve on Goppa's.
Xing's refinement applies to both constructions
and yields nonlinear codes over the same alphabet,
making possible a direct comparison in which
the codes of~\cite{E:stoc} come out ahead.

Like the codes of \cite{E:stoc} and (probably) \cite{Xing},
and unlike Goppa codes, our new codes $\CCC_D(h)$
are still very far from any practical use:
it is not even clear that one can efficiently
encode integers in $[1,(\#\CCC_D(h))^{1-o(1)}]$
or recognize whether a given word in $k^N$ is in that code,
let alone solve the error-correcting problem.

The rest of this paper is organized as follows.
We first review the nonlinear algebraic-geometry codes 
$\CC_m$, $\sC_D(h)$ of~\cite{Xing,E:stoc}.
We then combine these two constructions, and show that
our new codes $\CCC_D(h)$ attain the claimed improvement
over Xing's codes.

This work is supported in part by the National Science Foundation
(grant DMS-0200687).
\vspace*{2ex}

{\large\bf 2. Two variations on a theme of Goppa}

{\bf 2.1 Xing: nonlinear codes from derivatives of sections.}
Let $C$\/ be a (projective, smooth, irreducible)
algebraic curve of genus~$g$ defined over~$k$,
with $N$\/ rational points $P_1,\ldots,P_N$.
For each $j=1,\ldots,N$, choose a local uniformizing parameter
$t_j$ at~$P_j$, that is, a rational function on~$C$\/
vanishing to order~$1$ at~$P_j$.
(Different choices of $t_j$ will yield isomorphic codes.)
For any rational function $f\in k(C)$,
let $f_j^{(r)}$ ($r=0,1,2,\ldots$) be the $t_j^r$ coefficient
of its expansion in powers of~$t_j$; that is, $f_j^{(0)} = f(P_j)$,
and $f_j^{(1)}, f_j^{(2)}, f_j^{(3)}, \ldots$ are chosen inductively
so that $t_j^{-m} \bigl( f - \sum_{0\leq r<m} f_j^{(r)} t_j^r \bigr)$
is regular at~$P_j$ for all~$m$.

Let $D$\/ be a divisor on~$C$.  As in~\cite{Xing},
we simplify the exposition by assuming that the support of~$D$\/
is disjoint from $C(k)=\{P_1,\ldots,P_N\}$.
(Little is lost by this assumption, because any divisor on~$C$\/
is linearly equivalent to one satisfying the disjointness condition,
and linearly equivalent divisors yield equivalent codes.)
Note that we do {\em not}\/ assume that $\deg(D)<N$, and indeed
we shall use divisors of degree considerably larger than~$N$.
We denote by $\LL(D)$ the vector space of global sections of~$D$.
For any distinct $f,f'\in\LL(D)$, the difference $f-f'$ has a total of
$\deg(D)$ zeros in $C(\kbar)$, {\em counted with multiplicity}.
Goppa's construction does not exploit multiplicity,
using only the corollary that $f(P_j)=f'(P_j)$ holds
for at most $\deg(D)$ values of~$j$.
Multiple zeros are the key to Xing's improvement.

For $r\geq 0$ define $\phi_r:\LL(D) \ra k^N$ by
\begin{equation}
\label{eq:phi_r}
\phi_r(f) = \bigl( f_1^{(r)}, f_2^{(r)},\ldots, f_N^{(r)} \bigr).
\end{equation}
Recall that the {\em Hamming distance} $d(\cdot,\cdot)$ on~$k^N$
is defined by
\begin{equation}
\label{eq:d}
d(c,c') :=
\# \bigl\{ j\in\{1,2,\ldots,N\} \mid c_j\0 \neq c_j' \bigr\}.
\end{equation}
For each positive integer~$m$, Xing defines a code $\CC_m$
as follows. 
For $r=0,1,\ldots,m-1$, fix positive real numbers $\sigma_r < (q-1)/q$,
which will be specified later to optimize~$\CC_m$.
Choose $c_r\in k^N$ to maximize the size of
\begin{equation}
\label{eq:M_m}
\MM_m := \bigl\{
  f \in \LL(D) \mid d(c_r,\phi_r(f)) \leq \sigma_r N
  \
  {\rm for\ each}\ r=0,1,\ldots,m-1
\bigr\},
\end{equation}
and define
\begin{equation}
\label{eq:C_m}
\CC_m := \phi_m(\MM_m).
\end{equation}
As $(c_0,\ldots,c_{m-1})$ varies over all $q^{mN}$ possible choices,
the average size of $\MM_m$ equals $q^{-mN} \#(\LL(D))$
times the product of the sizes of closed Hamming balls
of radius $\sigma_0 N$, $\sigma_1 N$, \ldots, $\sigma_{m-1} N$.
The maximal $\MM_m$ is therefore at least as large, so
\begin{equation}
\label{eq:|M_m|}
\# \MM_m > q^{(H-o(1)) N} \#(\LL(D)),
\end{equation}
where the negative number $H$\/ is given in terms of
the entropy function~(\ref{eq:entropy}) by
\begin{equation}
\label{eq:H}
H = \sum_{r=0}^{m-1} (H_q(\sigma_r)-1).
\end{equation}
Now let $f,f'$ be distinct functions in~$\MM_m$,
and for each $r\geq 0$ let $I_r$ be the set of $i\in\{1,\ldots,N\}$
such that $f_j^{(r)} = {f'_j}^{(r)}$.
If $r<m$ then $\#I_r \geq (1-2\sigma_r) N\!$,
since $N-\#I_r$ is the distance between the words $\phi_r(f),\phi_r(f')$
in the same Hamming ball of radius $\sigma_r N$.
On the other hand, the total number of zeros of $f-f'$
at $P_1,\ldots,P_n$, counted with multiplicity, is
\pagebreak
\begin{eqnarray}
\sum_{s=0}^\infty \#(\cap_{r=0}^s I_r)
& \geq &
\sum_{s=0}^m \# \Bigl( \bigcap_{r=0}^s I_r \Bigr)
\nonumber
\\
& \geq &
(m+1)N - \sum_{s=0}^m \left(\sum_{r=0}^s (N - \# I_r)\right)
\label{eq:Isum}
\\
& = &
(m+1)N - \sum_{r=0}^m (m+1-r) (N - \# I_r)
\nonumber
\\
& \geq &
\Bigl(m - 2\sum_{r=0}^{m-1} (m+1-r) \sigma_r \Bigr)N
\, + \, \# I_m.
\nonumber
\end{eqnarray}
Since this number may not exceed the degree of~$D$, we deduce that
\begin{equation}
\label{eq:Imbound}
\# I_m \leq
\deg(D) - \Bigl(m - 2\sum_{r=0}^{m-1} (m+1-r) \sigma_r \Bigr)N.
\end{equation}
Define $d_0$, then, by
\begin{equation}
\label{eq:d0}
d_0 := \Bigl(m+1 - 2\sum_{r=0}^{m-1} (m+1-r) \sigma_r \Bigr)N
\, - \, \deg(D).
\end{equation}
Then $\#I_m \leq N-d_0$.  That is, $d(\phi_m(f),\phi_m(f')) \geq d_0$.
Therefore $\CC_m$ has minimum distance at least~$d_0$.
Assume that $\deg(D)$ is small enough that $d_0 > 0$:
\begin{equation}
\label{eq:Dbound}
\deg(D) < \Bigl(m+1 - 2\sum_{r=0}^{m-1} (m+1-r) \sigma_r \Bigr)N.
\end{equation}
Then $\#I_m < N$, that is, $\phi_m(f) \neq \phi_m(f')$.
Therefore $\phi_m$ is injective and $\#\CC_m = \#\MM_m$.

It remains to optimize $\deg(D)$ and $\sigma_r$ given $\delta = d_0/N$.
By Riemann-Roch, $\LL(D)$ is a \hbox{$k$-vector} space of dimension
at least $\deg(D)-g+1$, with equality if $D>2g-2$
(which will be the case for all $D$\/ that we use).
Combining this with (\ref{eq:d0})
and the estimate (\ref{eq:|M_m|}) on $\#\MM_m$, we find
\begin{equation}
\label{eq:Xingbd}
\frac{\log_q \#\CC_m}{N} + \frac{d_0}{N}
> 1 - \frac{g}{N}
+ \sum_{r=0}^{m-1} \bigl(H_q(\sigma_r) - 2(m+1-r)\sigma_r \bigr)
- o(1).
\end{equation}
The left-hand side is a lower bound on $R+\delta$ for the code $\CC_m$.
If we take all $\sigma_r=0$, the right-hand side reduces to $1-(g/N)$,
and we recover the Goppa bound~(\ref{eq:Goppa}) --- which is to be
expected because in this case
$\CC_m$ is equivalent to the Goppa code $\LL(D-m\sum_{i=1}^N (P_j))$
(this is most easily seen if we also choose each $c_r=0$).
Since the derivative of $H_q(\sigma) - 2(m+1-r)$ at $\sigma=0$
is $+\infty$ for any $m,r$,
there must exist positive $\sigma_r$ that improve on Goppa.
By differentiating each term of the sum in~(\ref{eq:Xingbd}),
Xing computed that the optimal choice is
\begin{equation}
\label{eq:sigma}
\sigma_r = (q-1) / (q^{2(m+1-r)} + q - 1),
\end{equation}
corresponding to $\sigma_r / (1-\sigma_r) = q^{-2(m+1-r)}(q-1)$.
Using the equivalent formula (\ref{eq:entropy1}) for $H_q$
and taking $i=m+1-r$,
the \hbox{$m$-th} term of~(\ref{eq:Xingbd}) is seen to equal
$\log_q(1+q^{-2i}(q-1))$ at the optimal~$\sigma_r$,
whence the codes $\CC_m$ attain
\begin{equation}
\label{eq:Xing_m}
R + \delta > 1 - \frac{1}{q_0-1}
  + \sum_{i=2}^{m+1} \log_q \left( 1 + \frac{q-1}{q^{2i}} \right)
  - o(1).
\end{equation}
Taking $m\ra\infty$ and $q=q_0^2$ recovers (\ref{eq:Xing});
without the hypothesis $q=q_0^2$, Xing's construction still improves
on (\ref{eq:Goppa}) by adding $c_q$ to the lower bound on $R+\delta$.

{\bf 2.2 Codes over $\Pr^1(k)$ from rational functions and sections.}

The codes we introduced in~\cite{E:stoc} use rational functions
on~$C$\/ instead of the global sections that comprise Goppa's codes.
Let $D$\/ be a divisor {\em of degree zero} on~$C$.
For a nonnegative integer~$h$, we define $\sM_D(h)$
to be the set of rational sections of degree $\leq h$\/
of the line bundle~$L_D$ associated to~$D$.
That is, $\sM_D(h) \subset k(C)$ consists of the zero function
together with the rational functions~$f$\/ on~$C$\/
whose divisor $(f)$ is of the form $E-D$, where $E$\/ is a divisor
whose positive and negative parts each have degree at most~$h$.
For instance, $\sM_0(h)$ is the set of rational functions
of degree at most~$h$.  Here, as opposed to \cite{E:stoc},
we do {\em not}\/ assume that $h<N/2$,
and indeed we shall use considerably larger~$h$.

To use these $\sM_D(h)$ for coding,
we need an upper bound on the number of solutions, with multiplicity,
of $f(P)=f'(P)$ for distinct rational sections $f,f'$ of~$L_D$,
given the degrees of~$f$ and~$f'$.
We define this multiplicity as follows.
For each $P\in C(\kbar)$ choose a rational function $\vph_P$
whose divisor has the same order at~$P$\/ as~$D$.
(The definition will be clearly independent
of the choice of~$\vph_P$.)
Then $P$\/ is a solution of $f=f'$
if the rational functions $\vph_P f$\/ and $\vph_P f'$
have the same value, finite or infinite, at~$P$.
In the former case, the multiplicity
is the valuation of $(\vph_P f) - (\vph_P f')$ at~$P$.
In the latter case, the multiplicity
is the valuation of $(\vph_P f)^{-1} - (\vph_P f')^{-1}$ at~$P$.
Of course if $P$\/ is not a solution of $f=f'$
then its multiplicity is~zero.

{\bf Proposition.}
{\em
Suppose $f,f'$ are distinct rational sections of~$L_D$,
with degrees $h,h'$.  For $P\in C(\kbar)$ let $m(P)$
be the multiplicity of~$P$\/ as a solution of $f=f'$.
Then $\sum_{P\in C(\kbar)} m(P) = h+h'$.
}

{\em Proof.}  For $P\in C(\kbar)$, if $P$\/ is a pole of $\vph_P f$,
let $\mu(P)$ be the multiplicity of this pole,
and otherwise set $\mu(P)=0$;
define $\mu'(P)$ likewise using $\vph_P f'$.
Then $h+h' = \sum_{P\in C(\kbar)} (\mu(P)+\mu'(P))$.
But we claim that $m(P) - (\mu(P)+\mu'(P))$
is the valuation at~$P$\/ of $\vph_P f - \vph_P f'$,
considered also as a rational section of~$L_D$.
This claim is immediate if neither $\vph_P f$\/ nor $\vph_P f'$
has a pole at~$P$\/; if just one of them has a pole there
then $\vph_P f - \vph_P f'$ has a pole of the same order,
which equals $\mu(P)+\mu'(P)$, while $m(P)=0$; 
finally, if both $\vph_P f, \vph_P f'$ have poles at~$P$,
then $m(P)>0$, and the claim follows from the identity
$$
(\vph_P f)^{-1} - (\vph_P f')^{-1} =
\frac{(\vph_P f') - (\vph_P f)} {(\vph_P f') (\vph_P f)}
$$
by taking valuations at~$P$.  This establishes our claim in all cases.
But the sum over $P\in C(\kbar)$ of $v_P\0(\vph_P f - \vph_P f')$
vanishes, since $\sum_{P\in C(\kbar)} v_P\0(f-f') = 0$ while
$\sum_{P\in C(\kbar)} v_P\0(\vph_P) = \deg(D)$
was also assumed to equal~zero.
Since
$$
\sum_P v_P\0(\vph_P f - \vph_P f')
=
\sum_P \bigl( m(P) - (\mu(P)+\mu'(P)) \bigr)
=
\sum_P m(P) - (h+h'),
$$
this proves that $m(P)=h+h'$.~~$\Box$

{\em Remark.}  This result is the analogue for rational sections
of the fact that a nonzero element of $\LL(D)$ has $\deg(D)$ zeros
counted with multiplicity.  It refines Proposition~1 of~\cite{E:stoc},
where we showed only that $f=f'$ has at most $h+h'$ solutions
not counting multiplicity.

In particular, if $2h<N$\/ then the evaluation map
$\phi_0: \sM_D(h) \ra \bigl(\Pr^1(k)\bigr)^{\!N}$ taking $f$\/ to
$\bigl( (\vph_{P_1} f) (P_1), \ldots, (\vph_{P_N} f) (P_N) \bigr)$
is injective, and its image is a code of length~$N$\/ over $\Pr^1(k)$
with minimal distance at least $N-2h$.  We call this code $\sC_D(h)$.
Note that it is $2h$, rather than $h$, that plays the role analogous
to the degree of the divisor on Goppa's construction; the notation $h$\/
should suggest both the {\em h}\/eight of a rational section of~$L_D$
and {\em h}\/alf of the degree of Goppa's divisor.

We also need the size of this code.  It turns out to be easier,
though still far from trivial, to estimate
not individual $\#\sM_D(h)$ but the average of $\#\sM_D(h)$
as $D$\/ ranges over (representatives of) the Jacobian $J_C$,
which is the group of equivalence classes of degree-zero divisors on~$C$.
We quote the following from \cite[Thm.1]{E:stoc}:\footnote{
  For individual codes, see \cite[Thm.2]{E:stoc},
  which gives the same estimate but only
  with a much higher threshold on~$\rho$.
  Since we later use Xing's technique,
  which requires an averaging argument, we would gain little
  by citing here an estimate on $\#\sM_D(h)$ free of averaging.
  Therefore we do not invoke \cite[Thm.2]{E:stoc}
  in our present application.
  }
\begin{quote}
{\em
If $C$\/ is asymptotically optimal (i.e., if $C$\/ varies
in a family of curves of genus \hbox{$g\ra\infty$}
with $N\sim(q^{1/2}-1)g$ rational points),
and for each~$C$\/ we choose $h$ with $\rho = \inf(h/N) > q/(q^2-1)$,
then the average over $D\in J_C$ of $\#\sM_D(h)$ is
\begin{equation}
\label{eq:Masymp}
\left(\frac{q+1}{q}\right)^{N \pm o_\rho(N)} q^{2h-g}.
\end{equation}
}
\end{quote}

Hence there exist $D$\/ for which $\sM_D(h)$
has size at least (\ref{eq:Masymp}),
which exceeds by a factor of $\bigl((q+1)/q\bigr)^{N-o_\rho(N)}$
the Riemann-Roch estimate on the size of $\LL(D)$ when $\deg(D)=2h$.
Note that we do not require that $2h<N$.
We shall use the result also for some $h \geq N/2$, 
in which case we cannot deduce a lower bound on~$\sC_D(h)$
(and anyway we have no nontrivial lower bound
on the minimal distance of~$\sC_D(h)$),
but will be able to construct another code $\CCC_D(h)$
by adapting Xing's technique.

\pagebreak

{\large\bf 3. New nonlinear codes over~$k$}

We again regard the evaluation map $\phi_0$ on $\sM_D(h)$
as the first of an infinite series $\phi_0,\phi_1,\phi_2,\ldots$
that records not just the values but also the derivatives
of rational sections of~$D$\/ at $P_1,\ldots,P_N$.
While $\phi_0$ takes values in $\bigl(\Pr^1(k)\bigr)^{\!N}$,
the $\phi_r$ for $r>0$ take values in $k^N$; this is why
we construct codes with alphabet~$k$\/ rather than $\Pr^1(k)$.
We define $\phi_1,\phi_2,\ldots$ as follows.
If $(\vph_{P_j} f) (P_j)=\infty$
then the \hbox{$j$-th} coordinate of $\phi_r(f)$
is the $t_j^r$ coefficient of the expansion of $(\vph_{P_j} f)^{-1}$
in powers of~$t_j$.  Otherwise that coordinate
is the $t_j^r$ coefficient of the expansion of $\vph_{P_j} f$.
In either case, $P_j$ is a solution of $f=f'$ of multiplicity~$m$
if and only if the \hbox{$j$-th} coordinates of $\phi_r(f), \phi_r(f')$
coincide for $0 \leq r < m$ but not for $r=m$.
As in~\cite{Xing}, we could have dispensed with the $\vph_{P_j}$
by choosing a linearly equivalent $D$\/ whose support is disjoint
from $\{P_1,\ldots,P_N\}$.

Fix positive $\sigma_0 < q/(q+1)$.
Choose $c_0 \in \bigl(\Pr^1(k)\bigr)^{\!N}$
to maximize the size of
\begin{equation}
\label{eq:MMM_m}
\MMM_D(h) := \bigl\{
  f \in \sM_D(h) \mid d(c_0,\phi_0(f)) \leq \sigma_0 N
\bigr\},
\end{equation}
and define
\begin{equation}
\label{eq:CCC}
\CCC_D(h) := \phi_1(\MMM_D(h)).
\end{equation}
As $c_0$ varies over all $(q+1)^N$ possible choices,
the average size of $\MMM_D(h)$ equals $(q+1)^{-N} \#(\MMM_D(h))$
times size of the closed Hamming ball of radius $\sigma_0 N$
in $\bigl(\Pr^1(k)\bigr)^{\!N}$.
The maximal $\MMM_D(h)$ is therefore at least as large, so
\begin{equation}
\label{eq:|MMM|}
\# \MMM_D(h) > (q+1)^{(H_{q+1}(\sigma_0)-1-o(1)) N} \#(\sM_D(h)).
\end{equation}
If moreover $C$\/ is asymptotically optimal
and $h \geq \rho N$\/ for some $\rho > q/(q^2-1)$,
then there exists~$D$\/ such that
$\#(\sM_D(h))$ is bounded below by~(\ref{eq:Masymp}).
We then have
\begin{equation}
\label{eq:|MMM|+}
\# \MMM_D(h) > q^{2h-g}
\exp \bigl[ N \bigl(
  \log(q+1) H_{q+1}(\sigma_0) - \log q - o_\rho(1)
\bigr) \bigr].
\end{equation}
We next fix $d_0>0$ and show that if
\begin{equation}
\label{eq:2h}
2h \leq (2-4\sigma_0)N - d_0
\end{equation}
then $\phi_1$ is injective on $\MMM_D(h)$
and its image $\CCC_D(h)$ is a code of minimum distance at least~$d_0$.
We have seen that for any distinct $f,f'\in\sM_D(h)$
the total multiplicity of solutions of $f=f'$ is at most~$2h$.  
If $f,f'\in\MMM_D(h)$ then $f(P_j)=f'(P_j)$
for at least $N - 2\sigma_0 N$ values of~$j$.
For at least
$$
N - 2\sigma_0 N - d(\phi_1(f),\phi_1(f'))
$$
of these,
the \hbox{$j$-th} coordinates of $\phi_1(f)$ and $\phi_1(f')$
also coincide, so $P_j$ is a solution of $f=f'$
of multiplicity at least~$2$.
Hence
$$
(2-4\sigma_0)N - d(\phi_1(f),\phi_1(f')) \leq 2h.
$$
Therefore if (\ref{eq:2h}) holds then $d(\phi_1(f),\phi_1(f')) \geq d_0$,
which proves are claim that $\phi_1$ maps $\MMM_D(h)$ injectively
to a code of minimum distance $\geq d_0$.
(This is of course a direct adaptation
of the case $m=1$ of Xing's argument,
which is Theorem~1.2 and \S II of his paper~\cite{Xing}.)

Finally we take $d_0 = \delta N$, let $2h = (2-4\sigma_0-\delta-o(1))N$,
and optimize $\sigma_0$.  Combining our results thus far, we have
\begin{equation}
\label{eq:CCCbd}
\frac{\log_q \#\CCC_m}{N} + \delta
> 1 - \frac{g}{N}
+ \log_q(q+1) \cdot H_{q+1}(\sigma_0) - 4\sigma_0
- o_\rho(1).
\end{equation}
If we took $\sigma_0=0$,
we would again recover the Goppa bound (\ref{eq:Goppamax}).\footnote{
  In \cite{E:stoc} we already noted, in a different way,
  that we could use $\sC_d(h)$ to construct nonlinear codes
  that are asymptotically as good as Goppa's, and cited this
  observation in one of our indirect comparisons
  between $\sC_d(h)$ and Goppa codes.
  }
Since the derivative with respect to~$\sigma_0$
of the bound~(\ref{eq:CCCbd}) is $+\infty$,
the optimal~$\sigma_0$ must improve on (\ref{eq:Goppamax})
as was the case for Xing's construction.
The resulting bound improves on Xing's because (\ref{eq:CCCbd})
involves the entropy function for an alphabet of $q+1$ letters
rather than~$q$.  Here we calculate that the optimal~$\sigma_0$ is
$1/(q^3+1)$, corresponding to $\sigma_0 / (1-\sigma_0) = q^{-3}$.
Substituting this into (\ref{eq:CCCbd}) we obtain
\begin{equation}
\label{eq:CCCopt}
\frac{\log_q \#\CCC_m}{N} + \delta
> 1 - \frac{g}{N}
+ \log_q \left( 1 + \frac1{q^3} \right)
- o_\rho(1).
\end{equation}
With this choice of $\sigma_0$, the ratio $2h/N$ easily exceeds
the threshold of $2q/(q^2-1)$ even for $\delta=1$ as long as $q_0\geq4$.
(We must in any event exclude $q=2$ or~$3$, because for those~$q$
none of the methods described here could yield
codes with positive $R,\delta$
even if asymptotically optimal curves were known.)

We have therefore proved:

{\bf Theorem.}
{\em
Let $q_0$ be a prime power, and $k$ a finite field of~$q^2$ elements.
For all positive $R_0,\delta_0$ satisfying (\ref{eq:feasible}),
and any $N_0$, there exist a curve $C/k$,
a degree-zero divisor $D$ on~$C$,
and a positive integer~$h$, such that $\CCC_D(h)$
is a code whose length $N$, transmission rate~$R$,
and error-detection rate $\delta=d/N$\/ satisfy
$N > N_0$, $R > R_0$, and $\delta > \delta_0$.
In particular, all positive $(R_0,\delta_0)$
satisfying (\ref{eq:feasible}) are asymptotically feasible.
}

As noted already, this construction corresponds to the case $m=1$
of Xing's codes.  One can formulate such a construction for any~$m$,
but with worse results, because only $\sigma_0$ increases,
and the resulting improvement falls off as $q^{-2m-2}/\log(q)$
when $m$ grows.  Can multiplicities of order $>2$ be exploited
to yield even better codes?

\end{document}